\documentclass[a4paper,10pt]{article}
\usepackage{a4wide}
\usepackage{amsmath}
\usepackage{amssymb}
\usepackage{amsthm}
\usepackage{latexsym}
\usepackage{graphicx}
\usepackage[english]{babel}
\usepackage{makeidx}

\newtheorem{theorem}{Theorem}[section]
\newtheorem{lemma}[theorem]{Lemma}
\newtheorem{corollary}[theorem]{Corollary}

\theoremstyle{definition}

\theoremstyle{remark}

\begin{document}
\title{Intersections of base rings associated to transversal polymatroids}
\author{Alin \c{S}tefan\\
 "Petroleum and Gas" University of Ploie\c sti, Romania}
\date{}
\maketitle

\begin{abstract}
In this paper, we study when the intersection of base rings
associated to some transversal polymatroids from \cite {SA} is the
base ring of a transversal polymatroid. This intersection  is a
Gorenstein ring and we compute its $a$-invariant.
\end{abstract}

\begin{quotation}


\end{quotation}

\section{Introduction}

The discrete polymatroids and their base rings are studied recently
in many papers  (see \cite{HH}, \cite{HHV}, \cite{V1}, \cite{V2}).
It is important to give conditions when the base ring associated to
a transversal polymatroid is Gorenstein (see \cite{HH}). In
\cite{SA} we introduced a class of such base rings. In this paper we
note that an intersection of such base rings (introduced in
\cite{SA}) is Gorenstein and give necessary and sufficient
conditions for the intersection of two base rings from \cite{SA} to
be still a base ring of a transversal polymatroid. Also, we compute
the $a$-invariant of those base rings. The results presented were
discovered by extensive computer algebra experiments performed with
{\it{Normaliz}} \ \cite{BK}.

\section{Preliminaries}

Let $n\in \mathbb{N},$ \ $n \geq 3,$ $\sigma \in S_{n},$ \
$\sigma=(1,2,\ldots, n)$ \ the cycle of length $n,$ \ $[n]:=\{1,
2,\ldots, n\},$ $\sigma^{t}[i]:=\{\sigma^{t}(1),\ldots,\sigma^{t}(i)\}$
for any $1\leq i \leq n-1$ \ and $\{e_{i}\}_{1 \leq i \leq n}$ \ be the
canonical base of $\mathbb{R}^{n}.$  For a vector $x\in
\mathbb{R}^{n}$, $x=(x_{1},\ldots,x_{n})$,  we will denote  $\mid x \mid
:= x_{1}+ \ldots + x_{n}.$ If $x^{a}$ is a monomial in $K[x_{1},
\ldots, x_{n}]$ we set $\log (x^{a})=a$. Given a set $A$ of monomials,
the $log \ set \ of \ A,$ denoted $\log (A),$ consists of all
$\log (x^{a})$ with $x^{a}\in A.$\\
If $A_{i}$ are some nonempty subsets of $[n]$ for $1\leq i\leq
m$, $m\geq 2$, ${\bf{\mathcal{A}}}=\{A_{1},\ldots,A_{m}\}$, then the
set of the vectors $\sum_{k=1}^{m} e_{i_{k}}$ with $i_{k} \in
A_{k}$ is the base of a polymatroid, called the \emph{transversal
polymatroid presented} by ${\bf{\mathcal{A}}}.$ The \emph{base ring} of
a transversal polymatroid presented by ${\bf{\mathcal{A}}}$
is the ring
\[K[{\bf{\mathcal{A}}}]:=K[x_{i_{1}}\cdot\cdot\cdot
x_{i_{m}}\ | \  i_{j}\in A_{j},1\leq j\leq m].\]
We know by \cite{HH} that the $K-$ algebra $K[{\bf{\mathcal{A}}}]$ is
normal and hence Cohen-Macaulay.
From \cite{SA} we know that the transversal polymatroid presented
by \[{\bf{\mathcal{A}}}=\{ A_{k} \ | \ A_{\sigma^{t}{(k)}}=[n],\ if
\ k\in [i]\cup \{n\} \ and \ A_{\sigma^{t}{(k)}}=[n]\setminus
\sigma^{t}[i], \ if \ k\in [n-1]\setminus [i]\}\] has the base ring
associated  $K[{\bf{\mathcal{A}}}]$ a Gorenstein ring for any
$1\leq i \leq n-1$ and $0\leq t \leq n-1$.\\
We put
\[\nu_{\sigma^{t}[i]}=-(n-i-1)\sum_{k=1}^{i}e_{\sigma^{t}(k)}+
(i+1)\sum_{k=i+1}^{n}e_{\sigma^{t}(k)}\]
for any $0\leq t \leq n-1$ and $1\leq i \leq n-1.$\\
The cone generated by the exponent vectors of the monomials defining
the base ring $K[{\bf{\mathcal{A}}}],$\\
$A:=\{ log(x_{j_{1}}\cdot\cdot\cdot x_{j_{n}}) \ | \
j_{k}\in A_{k}, k \in [n]\}\subset \mathbb{N}^{n},$
has the irreducible representation:
\[{\mathbb{R_{+}}}A= \bigcap_{a\in N}H^{+}_{a},\] where
$N=\{\nu_{\sigma^{t}[i]},\ \nu_{\sigma^{k}[n-1]} \ | \ 0\leq k \leq n-1\}.$ \\
It is easy to see that for any $1\leq i \leq n-1$ and $0\leq t \leq n-1$
\[A=\{\alpha \in  \mathbb{N}^{n} \ | \ \mid \alpha \mid =n,
\ 0\leq \alpha_{t+1}+\ldots + \alpha_{t+i}\leq i+1 \}, \ if \ i+t \leq n
\ \ \ \ \ \]
and \[A=\{\alpha \in  \mathbb{N}^{n} \ | \ \mid \alpha \mid =n, \
0\leq \sum_{s=1}^{i+t-n}\alpha_{s}+ \sum_{s=t+1}^{n}\alpha_{s}\leq i+1 \},
\ if \ i+t \ >  \ n .\]
We denote by $\{x_{i_{1}},\ldots,x_{i_{s}}\}^{r}$ all monomials of degree
\emph{r} with the indeterminates $x_{i_{1}},\ldots,x_{i_{s}}.$

\section{Intersection of cones of dimension $n$ with $n+1$ facets}
Let $r\geq 2$, $1\leq i_{1}, \ldots, i_{r} \leq n-2,$ $0=t_{1}\leq t_{2},
\ldots, t_{r}\leq n-1$ and we consider \emph{r} presentations of transversal
polymatroids:
\[{\bf{\mathcal{A}}}_{s}=\{A_{s, k} \ | \ A_{s, \sigma^{t_{2}}(k)}=[n],\ if \
k\in[i_{2}]\cup \{n\},\ A_{s, \sigma^{t_{2}}(k)}=[n]\setminus
\sigma^{t_{2}}[i_{2}],\ if \ k\in[n-1]\setminus [i_{2}]\}\]
for any $1\leq s \leq r.$
From \cite{SA} we know that the base rings $K[{\bf{\mathcal{A}}}_{s}]$ are
Gorenstein rings and the cones generated by the exponent vectors of the
monomials defining the base ring associated to a transversal polymatroid
presented by ${\bf{\mathcal{A}}}_{s}$ are :
\[{\mathbb{R_{+}}}A_{s}= \bigcap_{a\in N_{s}}H^{+}_{a},\]
\[where \ N_{s}=\{\nu_{\sigma^{t_{s}}[i_{s}]},\ \nu_{\sigma^
{k}[n-1]} \ | \ 0\leq k \leq n-1\}, \ A_{s}=\{ log(x_{j_{1}}\cdot\cdot\cdot
x_{j_{n}}) \ | \ j_{k}\in A_{s, k}, 1\leq k \leq n\}\subset \mathbb{N}^{n} \]
for any  $1\leq s \leq r.$\\
We denote by $K[A_{1}\cap \ldots \cap A_{r}],$ $K- algebra$ generated
by $x^{\alpha}$ with $\alpha \in A_{1}\cap \ldots \cap A_{r}.$\\
It is clear that the cone
\[{\mathbb{R_{+}}}(A_{1}\cap \ldots \cap A_{r})\subseteq{\mathbb{R_{+}}}A_{1}
\cap \ldots \cap {\mathbb{R_{+}}}A_{r}= \bigcap_{a\in N_{1}\cup \ldots \cup
N_{r}}H^{+}_{a}.\]
Conversely, since \[A_{1}\cap \ldots \cap A_{r}=\{\alpha \in \mathbb{N}^{n}
\ | \ \mid \alpha \mid=n, \ H_{\nu_{\sigma^{t_{s}}[i_{s}]}}(\alpha)\geq 0,
\ for \ any \ 1\leq s \leq r\},\]
we have that \[{\mathbb{R_{+}}}(A_{1}\cap \ldots \cap A_{r})\supseteq
\bigcap_{a\in N_{1}\cup \ldots \cup N_{r}}H^{+}_{a}\]
and so,  \[{\mathbb{R_{+}}}(A_{1}\cap \ldots \cap A_{r})=
\bigcap_{a\in N_{1}\cup \ldots \cup N_{r}}H^{+}_{a}.\]
We claim that the intersection \[\bigcap_{a\in N_{1}\cup \ldots \cup
N_{r}}H^{+}_{a}\] is the irreducible representation of the cone
${\mathbb{R_{+}}}(A_{1}\cap \ldots \cap A_{r}).$\\
We prove by induction on $r\geq 1.$\ If $r=1$, then the intersection
\[\bigcap_{a\in N_{1}}H^{+}_{a}\] is the irreducible representation of the
cone ${\mathbb{R_{+}}}(A_{1}).$ (see \cite{SA}, Lemma \ 4.1)\\ If $r>1$ then
we have two cases to study:\\
$1)$ if we delete, for some $0\leq k \leq n-1,$ \ the hyperplane with
the normal $\nu_{\sigma^{k}[n-1]}$, then a coordinate
of a $\log (x_{j_{1}}\cdot\cdot\cdot
 x_{j_{i_{1}}}x_{j_{i_{1}+1}}\cdot\cdot\cdot x_{j_{n-1}}x_{j_{n}})$
would be negative, which is impossible;\\
$2)$ if we delete, for some $1\leq s \leq r$,  the hyperplane with
the normal $\nu_{\sigma^{t_{s}}[i_{s}]},$ then \[\bigcap_{a\in
N_{1}\cup \ldots \cup N_{s-1} \cup N_{s+1}\ldots N_{r}}H^{+}_{a}\]
is by induction the irreducible representation
of the cone ${\mathbb{R_{+}}}(A_{1}\cap \ldots \cap A_{s-1}\cap
A_{s+1}\ldots \cap A_{r})$, which is diffrent from ${\mathbb{R_{+}}}
(A_{1}\cap \ldots \cap A_{r}).$
Hence, the intersection \[\bigcap_{a\in N_{1}\cup \ldots \cup
N_{r}}H^{+}_{a}\] is the irreducible representation of the cone
${\mathbb{R_{+}}}(A_{1}\cap \ldots \cap A_{r}).$
\begin{lemma}
The $K-$ algebra $K[A_{1}\cap \ldots \cap A_{r}]$ is a Gorenstein ring.
\end{lemma}
\begin{proof}
We will show that the canonical module $\omega_{K[A_{1}\cap \ldots
\cap A_{r}]}$ is generated by $(x_{1}\cdot\cdot\cdot x_{n})K[A_{1}
\cap \ldots \cap A_{r}].$
Since the semigroups $\mathbb{N}(A_{t})$ are  normal for any $1\leq
t \leq r$, it follows that $\mathbb{N}(A_{1}\cap \ldots \cap A_{r})$
is normal. Then the $K-$ algebra $K[A_{1}\cap \ldots \cap A_{r}]$
is normal \  (see \cite{BH} Theorem 6.1.4. p. 260 \ ) \ and using
the $Danilov-Stanley$ theorem we get that the canonical
module $\omega_{K[A_{1}\cap \ldots \cap A_{r}]}$ is \[\omega_
{K[A_{1}\cap \ldots \cap A_{r}]}=(\{x^{\alpha}
\ | \ \alpha \in \mathbb{N}(A_{1}\cap \ldots \cap A_{r}) \cap
ri(\mathbb{R}_{+}(A_{1}\cap \ldots \cap A_{r}))\}).\]
Let $d_{t}$ be the greatest common divisor of $n \ and \
i_{t}+1, \ gcd(n, i_{t}+1)=d_{t},$ for any $1\leq t \leq r.$\\
For any $1\leq s \leq r$, there exist two possibilities for
the equation of the facet $H_{\nu_{\sigma^{t_{s}}[i_{s}]}}$ :\\
$1)$ If $i_{s}+t_{s}\leq n,$ then the equation of the facet
$H_{\nu_{\sigma^{t_{s}}[i_{s}]}}$ is:
\[H_{\nu_{\sigma^{t_{s}}[i_{s}]}}(y): \ \frac{(i_{s}+1)}{d_{s}}
\sum_{k=1}^{t_{s}}y_{k}-\frac{(n-i_{s}-1)}{d_{s}}\sum_{k=t_{s}+1}^
{t_{s}+i_{s}}y_{k}+\frac{(i_{s}+1)}{d_{s}}\sum_{k=t_{s}+i_{s}+1}^{n}y_{k}
=0. \ \ \ \ \ \ \ \ \ \ \ \ \ \ \ \ \ \]
$2)$ If $i_{s}+t_{s}\ > \ n,$ then the equation of the facet
$H_{\nu_{\sigma^{t_{s}}[i_{s}]}}$ is:
\[H_{\nu_{\sigma^{t_{s}}[i_{s}]}}(y): \ -\frac{(n-i_{s}-1)}{d_{s}}
\sum_{k=1}^{i_{s}+t_{s}-n}y_{k}+ \frac{(i_{s}+1)}{d_{s}}\sum_{k=i_{s}+t_{s}
-n+1}^{t_{s}}y_{k}-\frac{(n-i_{s}-1)}{d_{s}}\sum_{k=t_{s}+1}^{n}y_{k} =0.\]
The  relative interior of the cone $\mathbb{R}_{+}(A_{1}\cap \ldots \cap A_{r})$
is:  \[ri(\mathbb{R}_{+}(A_{1}\cap \ldots \cap A_{r}))=\{y \in \mathbb{R}^{n}
\ | \ y_{k} \ > \ 0, \ H_{\nu_{\sigma^{t_{s}}[i_{s}]}}(y) \ > \ 0
\ for \ any \ 1\leq k \leq n \ and \ 1\leq s \leq r\}\]
We will show that \[\mathbb{N}(A_{1}\cap \ldots \cap A_{r}) \cap
ri(\mathbb{R}_{+}(A_{1}\cap \ldots \cap A_{r}))= \ (1,\ldots, 1) \ + \
(\mathbb{N}(A_{1}\cap \ldots \cap A_{r}) \cap \mathbb{R}_{+}(A_{1}\cap
\ldots \cap A_{r})).\]
It is clear that $ri(\mathbb{R}_{+}(A_{1}\cap \ldots \cap A_{r})) \
\supset \ (1,\ldots, 1) \ + \mathbb{R}_{+}(A_{1}\cap \ldots \cap A_{r}).$\\
If $(\alpha_{1}, \alpha_{2}, \ldots, \alpha_{n})\in \mathbb{N}
(A_{1}\cap \ldots \cap A_{r}) \cap ri(\mathbb{R}_{+}(A_{1}\cap \ldots \cap A_{r})),$
then $\alpha_{k} \geq 1$ for any $1\leq k\leq n$ and for any $1\leq s \leq r$ we
have \[\frac{(i_{s}+1)}{d_{s}}\sum_{k=1}^{t_{s}}\alpha_{k}- \frac{(n-i_{s}-1)}
{d_{s}}\sum_{k=t_{s}+1}^{t_{s}+i_{s}}\alpha_{k}+\frac{(i_{s}+1)}
{d_{s}}\sum_{k=t_{s}+i_{s}+1}^{n}\alpha_{k} \geq 1, \ if \ i_{s}+t_{s}\leq n
\ \ \ \ \ \ \ \ \ \ \ \ \ \ \ \ \ \] or
\[-\frac{(n-i_{s}-1)}{d_{s}}\sum_{k=1}^{i_{s}+t_{s}-n}\alpha_{k}+ \frac{(i_{s}+1)}
{d_{s}}\sum_{k=i_{s}+t_{s}-n+1}^{t_{s}}\alpha_{k}-\frac{(n-i_{s}-1)}
{d_{s}}\sum_{k=t_{s}+1}^{n}\alpha_{k} \geq 1, \ if \ i_{s}+t_{s}\ > \ n\] and
\[ \sum_{k=1}^{n}\alpha_{k}=t \ n \ for \ some \ t\geq 1.\]
We claim that there exist $(\beta_{1}, \beta_{2}, \ldots, \beta_{n})
\in \mathbb{N}(A_{1}\cap \ldots \cap A_{r}) \cap \mathbb{R}_{+}
(A_{1}\cap \ldots \cap A_{r})$ such that $(\alpha_{1}, \alpha_{2}, \ldots, \alpha_{n})
=(\beta_{1}+1, \beta_{2}+1, \ldots, \beta_{n}+1).$
Let $\beta_{k}=\alpha_{k}-1$ for all $1\leq k \leq n.$\\
It is clear that $\beta_{k}\geq 0$ and for any $1\leq s \leq r,$
\[H_{\nu_{\sigma^{t_{s}}[i_{s}]}}(\beta)=H_{\nu_{\sigma^{t_{s}}[i_{s}]}}(\alpha)-
H_{\nu_{\sigma^{t_{s}}[i_{s}]}}(1,\ldots,1)=H_{\nu_{\sigma^{t_{s}}[i_{s}]}}
(\alpha)-\frac{n}{d_{s}}.\]
If $H_{\nu_{\sigma^{t_{s}}[i_{s}]}}(\beta)=j_{s}, \ for \ some \ 1\leq s \leq r \ and
 \ 1\leq j_{s} \leq \frac{n}{d_{s}}-1,$  then we will get a contadiction.
Indeed, since $n$ divides $\sum_{k=1}^{n}\alpha_{k},$ it follows $\frac{n}{d_{s}}$
divides $j_{s},$  which is false.\\ Hence, we have
$(\beta_{1}, \beta_{2}, \ldots, \beta_{n})\in \mathbb{N}(A_{1}\cap \ldots \cap A_{r})
\cap \mathbb{R}_{+}(A_{1}\cap \ldots \cap A_{r})$ and $(\alpha_{1}, \alpha_{2},
\ldots, \alpha_{n})\in \mathbb{N}(A_{1}\cap \ldots \cap A_{r}) \cap ri
(\mathbb{R}_{+}(A_{1}\cap \ldots \cap A_{r}))$.\\ Since $\mathbb{N}(A_{1}\cap \ldots
\cap A_{r}) \cap ri(\mathbb{R}_{+}(A_{1}\cap \ldots \cap A_{r}))= \ (1,\ldots, 1) \
+ \ (\mathbb{N}(A_{1}\cap \ldots \cap A_{r}) \cap \mathbb{R}_{+}(A_{1}\cap \ldots
\cap A_{r})),$ we get that $\omega_{K[A_{1}\cap \ldots \cap A_{r}]}= \
(x_{1}\cdot\cdot\cdot x_{n})K[A_{1}\cap \ldots \cap A_{r}].$
\end{proof}
Let $S$ be a standard graded $K-algebra$ over a field $K$. Recall that the
$a-invariant$ of $S,$ denoted $a(S),$ is the degree as a rational function of the
Hilbert series of $S,$ see for instance
$(\cite{V}, \ p. \ 99).$ If $S$ is $Cohen-Macaulay$ and $\omega_{S}$ is the
canonical module of $S,$ then \[a(S)= \ - \ min \ \{i\ | \ (\omega_{S})_{i}
\neq 0\},\] see (\cite{BH}, \ p. \ 141) and (\cite{V}, \ Proposition 4.2.3).
In our situation $S = K[A_{1}\cap \ldots \cap A_{r}]$ is normal  and
consequently $Cohen-Macaulay,$ thus this formula applies. As consequence of
$Lemma \ 3.1.$ we have the following:
\begin{corollary}
The $a-invariant$ of $K[A_{1}\cap \ldots \cap A_{r}]$ is $a(K[A_{1}\cap \ldots
\cap A_{r}])= \ -1.$
\end{corollary}
\begin{proof}
Let $\{x^{\alpha_{1}}, \ldots, x^{\alpha_{q}}\}$ be  the  generators  of
$K-algebra$ $K[A_{1}\cap \ldots \cap A_{r}].$ \  $K[A_{1}\cap \ldots
\cap A_{r}]$ is standard graded algebra with the grading
\[K[A_{1}\cap \ldots
\cap A_{r}]_{i}=\sum_{|c|=i}K(x^{\alpha_{1}})^{c_{1}}\cdot \cdot \cdot
(x^{\alpha_{q}})^{c_{q}}, \ where \ \mid c \mid =c_{1}+\ldots+c_{q}. \]
Since $\omega_{K[A_{1}\cap \ldots
\cap A_{r}]}= \ (x_{1}\cdot\cdot\cdot x_{n})K[A_{1}\cap \ldots
\cap A_{r}]$ it follows that $ min \ \{i\ | \ (\omega_{K[A_{1}\cap \ldots
\cap A_{r}]})_{i} \neq 0\}=1$ and so,  $a(K[A_{1}\cap \ldots
\cap A_{r}])= \ -1.$
\end{proof}

\section{When $K[A\cap B]$ is the base ring associated to some transversal polymatroid?}

Let $n\geq 2,$ and we consider two transversal polymatroids presented by:
${\bf{\mathcal{A}}}=\{A_{1},\ldots, A_{n}\}$ respectively ${\bf{\mathcal{B}}}=
\{B_{1},\ldots, B_{n}\}.$ Let $A$ and $B$ be the set of exponent vectors
of monomials defining the base rings $K[{\bf{\mathcal{A}}}],$ respectively
$K[{\bf{\mathcal{B}}}]$ and $K[A\cap B]$ the $K- algebra$ generated
by $x^{\alpha}$ with $\alpha \in A\cap B.$\\
{\bf{Question:}} There exists a transversal polymatroid such that
its base ring is the $K-algebra$ $K[A\cap B]$?\\
In the following we will give two suggestive examples.\\
{\bf{\emph{Example} $1$}.} \ Let $n=4,$ ${\bf{\mathcal{A}}}=\{A_{1},
A_{2}, A_{3}, A_{4}\},$ ${\bf{\mathcal{B}}}=\{B_{1}, B_{2}, B_{3},
B_{4}\},$ where $A_{1}=A_{4}=B_{2}= B_{3}=\{1, 2, 3, 4\},$
$A_{2}=A_{3}=\{2, 3, 4\},$  $B_{1}=B_{4}=\{1, 3, 4\}$ and
$K[{\bf{\mathcal{A}}}],$ $K[{\bf{\mathcal{B}}}]$ the base rings
associated to transversal polymatroids presented by
${\bf{\mathcal{A}}},$ respectively ${\bf{\mathcal{B}}}.$ It is easy
to see that the generators set of $K[{\bf{\mathcal{A}}}],$
respectively $K[{\bf{\mathcal{B}}}]$ are given by $A=\{y\in
\mathbb{N}^{4} \ | \  \mid y \mid =4, \ 0\leq y_{1}\leq 2, \
y_{k}\geq 0, \ 1\leq k \leq 4\},$ respectively $B=\{y\in
\mathbb{N}^{4} \ | \ \mid y \mid=4, \ 0\leq y_{2}\leq 2, \ y_{k}\geq
0, \ 1\leq k \leq 4\}.$ We show that the $K-$ algebra $K[A\cap B]$
is the base ring of the transversal polymatroid presented by
${\bf{\mathcal{C}}}=\{C_{1}, C_{2}, C_{3}, C_{4}\},$ where
$C_{1}=C_{4}=\{1, 3, 4\}, \ C_{2}=C_{3}=\{2, 3, 4\}.$\\
 Since the
base ring associated to the transversal polymatroid presented by
${\bf{\mathcal{C}}}$ has the exponent set $C=\{y\in \mathbb{N}^{4} \
| \ \mid y \mid=4, \ 0\leq y_{1}\leq 2, \ 0\leq y_{2}\leq 2, \
y_{k}\geq 0, \ 1\leq k \leq 4\},$ it follows that $K[A\cap
B]=K[{\bf{\mathcal{C}}}].$
Thus, in this example $K[A\cap B]$ is the base ring of a transversal polymatroid.\\
{\bf{\emph{Example} $2$}.} \ Let $n=4,$ ${\bf{\mathcal{A}}}=\{A_{1}, A_{2}, A_{3}, A_{4}\},$
${\bf{\mathcal{B}}}=\{B_{1}, B_{2}, B_{3}, B_{4}\}$ where
$A_{1}=A_{2}=A_{4}=B_{1}=B_{2}=B_{3}=\{1, 2, 3, 4\},$ $A_{3}=\{3, 4\},$
$B_{4}=\{1, 4\}$ and  $K[{\bf{\mathcal{A}}}],$ $K[{\bf{\mathcal{B}}}]$
the base rings associated to the transversal polymatroids presented by
${\bf{\mathcal{A}}},$ respectively ${\bf{\mathcal{B}}}.$
It is easy to see that the generators set of $K[{\bf{\mathcal{A}}}],$
respectively $K[{\bf{\mathcal{B}}}]$ are $A=\{y\in \mathbb{N}^{4} \ |
\ \mid y \mid=4, \ 0\leq y_{1}+y_{2}\leq 3, \ y_{k}\geq 0, \ 1\leq k \leq 4\},$
respectively $B=\{y\in \mathbb{N}^{4} \ | \ \mid y \mid=4, \
0\leq y_{2}+y_{3}\leq 3, \ y_{k}\geq 0, \ 1\leq k \leq 4\}.$
We claim that there exists no transversal polymatroid ${\bf{\mathcal{P}}}$
such that $K-$ algebra $K[A\cap B]$ is its base ring. Suppose,
on the contrary, let ${\bf{\mathcal{P}}}$ be presented by ${\bf{\mathcal{C}}}
=\{C_{1}, C_{2}, C_{3}, C_{4}\}$ with each $C_{k}\subset [4].$ Since $(3,0,1,0),
(3,0,0,1) \in {\bf{\mathcal{P}}}$ and $(3,1,0,0)\notin {\bf{\mathcal{P}}},$
we may assume that changing the numerotation of  $\{C_{i}\}_{i=\overline{1,4}}$
that  $1\in C_{1}, 1\in C_{2}, 1\in C_{4}$ and $C_{3}=\{3, 4\}.$ Since
$(0,3,0,1)\in {\bf{\mathcal{P}}},$ we assume that $2\in C_{1}, 2\in C_{2},
2\in C_{4}.$ Hence $(0,3,1,0)\in {\bf{\mathcal{P}}},$ a contradiction.\\

Let $1\leq i_{1}, i_{2} \leq n-2,$ $0\leq t_{2}\leq n-1$ and
$\tau \in S_{n-2},$ $\tau=(1,2,\ldots, n-2)$ \ the cycle of length $n-2.$
\ We consider two transversal polymatroids presented by:
\[{\bf{\mathcal{A}}}=\{A_{k} \ | \ A_{k}=[n],\ if \ k\in[i_{1}]\cup \{n\},
\ A_{k}=[n]\setminus [i_{1}],\ if \ k\in[n-1]\setminus [i_{1}]\}
\ \ \ \ \ \ \ \ \ \ \ \ \ \ \ \] and
\[{\bf{\mathcal{B}}}=\{B_{k} \ | \ B_{\sigma^{t_{2}}(k)}=[n],\ if \
k\in[i_{2}]\cup \{n\},\ B_{\sigma^{t_{2}}(k)}=[n]\setminus
\sigma^{t_{2}}[i_{2}],\ if \ k\in[n-1]\setminus [i_{2}]\}\] such that $A,$
respectively $B$ is the exponent vectors of the monomials defining the
base rings associated to transversal polymatroid presented by
${\bf{\mathcal{A}}},$ respectively ${\bf{\mathcal{B}}}.$
From \cite{SA} we know that the base rings $K[{\bf{\mathcal{A}}}],$
respectively $K[{\bf{\mathcal{B}}}]$ are Gorenstein rings and the cones
generated by the exponent vectors of the monomials defining the base ring
associated to the transversal polymatroids presented by ${\bf{\mathcal{A}}},$
respectively ${\bf{\mathcal{B}}}$ are :
\[{\mathbb{R_{+}}}A= \bigcap_{a\in N_{1}}H^{+}_{a}, \ \ \ \
{\mathbb{R_{+}}}B= \bigcap_{a\in N_{2}}H^{+}_{a},\] where $N_{1}=
\{\nu_{\sigma^{0}[i_{1}]},\ \nu_{\sigma^{k}[n-1]} \ | \ 0\leq k \leq
n-1\}, \ \ N_{2}=\{\nu_{\sigma^{t_{2}}[i_{2}]},\ \nu_{\sigma^{k}[n-1]}
\ | \ 0\leq k \leq
n-1\},$ \[A=\{ log(x_{j_{1}}\cdot\cdot\cdot x_{j_{n}}) \ | \
j_{k}\in A_{k}, 1\leq k \leq
n\}\subset \mathbb{N}^{n}\ and \ B=\{ log(x_{j_{1}}\cdot\cdot\cdot
x_{j_{n}}) \ | \ j_{k}\in B_{k}, 1\leq k \leq
n\}\subset \mathbb{N}^{n}.\]
It is easy to see that $A=\{\alpha \in  \mathbb{N}^{n} \ | \
\ 0\leq \alpha_{1}+\ldots + \alpha_{i_{1}}\leq i_{1}+1 \ and
\ \mid \alpha \mid=n \} $ and\\
$B=\{\alpha \in  \mathbb{N}^{n} \ | \ \ 0\leq
\alpha_{t_{2}+1}+\ldots + \alpha_{t_{2}+i_{2}}\leq i_{2}+1 \ and \
\mid \alpha \mid=n \}, \ if \ i_{2}+t_{2} \leq n $ \\ or\\
\[B=\{\alpha \in  \mathbb{N}^{n} \ | \ 0\leq
\sum_{s=1}^{i_{2}+t_{2}-n}\alpha_{s}+ \sum_{s=t_{2}+1}^{n}\alpha_{s}
\leq i_{2}+1 \ and \ \mid \alpha \mid=n \}, \ if \ i_{2}+t_{2} \geq n .
\ \ \ \ \ \ \ \ \ \ \ \ \ \ \ \ \ \ \ \ \]

For any base ring $K[{\bf{\mathcal{A}}}]$ of a transversal polymatroid
presented by ${\bf{\mathcal{A}}}=\{A_{1},\ldots ,A_{n}\}$ we \\
associate a $(n\times n)$ square tiled by closed unit subsquares,
called {\bf{boxes}}, colored with colors,
$"white"$ and $"black",$ as follows:
the box of coordinate $(i,j)$ is $"white"$ if $j\in A_{i},$
otherwise the box is $"black".$ We will call this square  the
{\bf{polymatroidal diagram}} associated to the presentation
${\bf{\mathcal{A}}}=\{A_{1},\ldots ,A_{n}\}.$\\

Next we give necesary and sufficient conditions such that the
$K-algebra$ $K[A\cap B]$ is the base ring associated to some
transversal polymatroid.

\begin{theorem}
Let $1\leq i_{1}, i_{2} \leq n-2,$ $0\leq t_{2}\leq n-1.$ We
consider two presentation of transversal poymatroids presented by:
${\bf{\mathcal{A}}}=\{A_{k} \ | \ A_{k}=[n],\ if \ k\in[i_{1}] \cup
\{n\},\ A_{k}=[n]\setminus [i_{1}],\ if \ k\in[n-1]\setminus
[i_{1}]\}$ and ${\bf{\mathcal{B}}}=\{B_{k} \ | \
B_{\sigma^{t_{2}}(k)}=[n],\ if \ k\in[i_{2}]\cup \{n\},\
B_{\sigma^{t_{2}}(k)}=[n]\setminus \sigma^{t_{2}}[i_{2}],\ if \
k\in[n-1]\setminus [i_{2}]\}$ such that $A,$ respectively $B$ are
the exponent vectors of the monomials defining the base ring
associated to the transversal polymatroids presented by
${\bf{\mathcal{A}}},$ respectively ${\bf{\mathcal{B}}}.$\\
Then, the $K-algebra$ $K[A\cap B]$ is the base ring associated to a
transversal polymatroid if and only if one of the
following conditions hold:\\
$a)$ \ \ $i_{1}=1.$ \\
$b)$ \ \ $i_{1}\geq 2$ and $t_{2}=0.$\\
$c)$ \ \ $i_{1}\geq 2$ and $t_{2}=i_{1}.$\\
$d)$ \ \ $i_{1}\geq 2,$  $1\leq t_{2} \leq i_{1}-1$ and
$i_{2}\in \{1,\ldots,i_{1}-t_{2}\}\cup \{n-t_{2},\ldots,n-2\};$\\
$e)$ \ \ $i_{1}\geq 2,$  $i_{1}+1\leq t_{2} \leq n-1$ and
$i_{2}\in \{1,\ldots,n-t_{2}\}
\cup \{n-t_{2}+i_{1},\ldots,n-2\}.$
\end{theorem}

The proof follows from the following three lemmas.

\begin{lemma}
Let $A$ and $B$ be like above. If $i_{1}\geq 2,$ $1\leq t_{2}
\leq i_{1}-1,$ then  the $K-algebra$ $K[A\cap B]$ is the base
ring associated to some transversal polymatroid if and only
if $i_{2}\in \{1,\ldots,i_{1}-t_{2}\}\cup \{n-t_{2},\ldots,n-2\}.$
\end{lemma}

\begin{proof}
$"\Leftarrow"$
Let $i_{2}\in \{1,\ldots,i_{1}-t_{2}\}\cup \{n-t_{2},\ldots,n-2\}.$
We will prove that there exists a transversal polymatroid ${\bf{\mathcal{P}}}$
presented by  ${\bf{\mathcal{C}}}=\{C_{1},\ldots,C_{n}\}$ such
that the base ring associated to ${\bf{\mathcal{P}}}$ is $K[A\cap B].$\\
We have two cases to study:\\
$\bf Case \ 1.$ If $i_{2}+t_{2}\leq i_{1},$ then let ${\bf{\mathcal{P}}}$
be the transversal polymatroid presented by  ${\bf{\mathcal{C}}}=
\{C_{1},\ldots,C_{n}\},$ where
\[C_{1}=\ldots=C_{i_{2}}=C_{n}=[n], \ \ \ \ \ \ \]
\[C_{i_{2}+1}=\ldots =C_{i_{1}}=[n]\setminus \sigma^{t_{2}}[i_{2}],\]
\[C_{i_{1}+1}=\ldots =C_{n-1}=[n]\setminus [i_{1}]. \ \]
The polymatroidal diagram associated is the following:

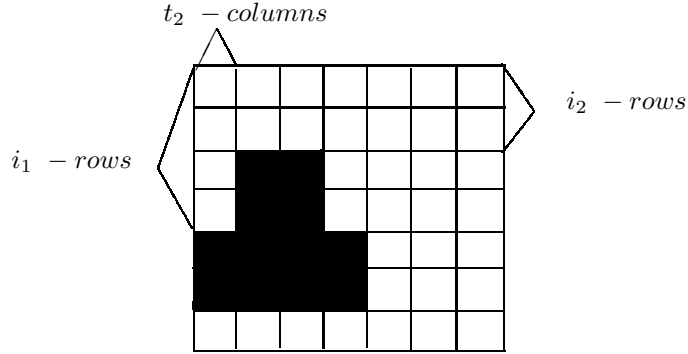
\unitlength 1mm 
\linethickness{0.4pt}
\ifx\plotpoint\undefined\newsavebox{\plotpoint}\fi 
\begin{picture}(0,67)(10,0)
\put(54,18.25){\framebox(40.75,37.75)[cc]{}}
\put(59.5,55.5){\line(0,-1){37}}
\put(65.25,56){\line(0,-1){37.25}}
\put(71,56){\line(0,-1){37.5}}
\put(76.75,55.75){\line(0,-1){37}}
\put(82.5,55.75){\line(0,-1){37.5}}
\put(88.5,56){\line(0,-1){37.75}}
\put(54,50.5){\line(1,0){41}}
\put(54.25,44.75){\line(1,0){40.5}}
\put(54,39.75){\line(1,0){40.75}}
\put(54,34){\line(1,0){40.75}}
\put(53.75,29.25){\line(1,0){41}}
\put(54,23.5){\line(1,0){40.5}}
\put(53.75,23.75){\rule{.5\unitlength}{1\unitlength}}
\put(54,23.75){\rule{22.75\unitlength}{10\unitlength}}
\put(59.5,34){\rule{11.5\unitlength}{10.5\unitlength}}
\multiput(54,56)(-.03368794,-.09574468){141}{\line(0,-1){.09574468}}
\multiput(49.25,42.5)(.03358209,-.05970149){134}{\line(0,-1){.05970149}}
\put(54.25,55.5){\line(1,2){2.75}}
\multiput(57,61)(.0333333,-.0633333){75}{\line(0,-1){.0633333}}
\multiput(70.75,45)(.03125,-.03125){8}{\line(0,-1){.03125}}
\multiput(94.5,56.25)(.03373016,-.04960317){126}{\line(0,-1){.04960317}}
\multiput(98.75,50)(-.03373016,-.04365079){126}{\line(0,-1){.04365079}}
\put(103,50.25){$i_{2} \ - rows$}
\put(30,42.5){$i_{1} \ -rows$}
\put(50,62){$t_{2} \ -columns$}
\end{picture}

It is easy to see that the base ring associated to the transversal
polymatroid ${\bf{\mathcal{P}}}$ presented by ${\bf{\mathcal{C}}}$
is generated by the following set of monomials:
\[\{x_{t_{2}+1},\ldots,x_{t_{2}+i_{2}}\}^{i_{2}+1-k}\{x_{1},\ldots,
x_{t_{2}},x_{t_{2}+i_{2}+1},\ldots,x_{i_{1}}\}^{i_{1}-i_{2}+k-s}
\{x_{i_{1}+1},\ldots,x_{n}\}^{n-1-i_{1}+ s}\] for any
$0\leq k \leq i_{2}+1$ and $0\leq s \leq i_{1}-i_{2}+k.$
If $x^{\alpha}\in K[{\bf{\mathcal{C}}}],$
$\alpha=(\alpha_{1},\ldots,\alpha_{n})\in {\mathbb{N}}^{n},$
then there exists $0\leq k \leq i_{2}+1$ and $0\leq s \leq i_{1}-i_{2}+k$
such that
\[\alpha_{t_{2}+1}+\ldots + \alpha_{t_{2}+i_{2}}= i_{2}+1-k \ and \
\alpha_{1}+\ldots + \alpha_{i_{1}}= i_{1}+1-s\]
and thus, $K[{\bf{\mathcal{C}}}]\subset K[A\cap B].$\\
Conversely, \ if $\alpha \in A\cap B$ then $\alpha_{t_{2}+1}+
\ldots + \alpha_{t_{2}+i_{2}}\leq i_{2}+1 \ and \ \alpha_{1}+\ldots +
\alpha_{i_{1}}\leq i_{1}+1$
and thus there exists $0\leq k \leq i_{2}+1$ and $0\leq s \leq i_{1}-i_{2}+k$
such that
\[x^{\alpha}\in \{x_{t_{2}+1},\ldots,x_{t_{2}+i_{2}}\}^{i_{2}+1-k}
\{x_{1},\ldots,x_{t_{2}}, x_{t_{2}+i_{2}+1},\ldots,x_{i_{1}}\}^
{i_{1}-i_{2}+k-s}\{x_{i_{1}+1},\ldots,x_{n}\}^{n-1-i_{1}+ s}.\]
Thus, $K[{\bf{\mathcal{C}}}]\supset K[A\cap B]$ and so $K[{\bf{\mathcal{C}}}]= K[A\cap B].$\\
${\bf Case} \ 2.$ If $i_{2}+t_{2}>i_{1},$ then it follows that $i_{2}\geq n-t_{2}$
and  $n-i_{2}\leq t_{2}\ < \ i_{1}.$
Let ${\bf{\mathcal{P}}}$ be the transversal polymatroid presented by
${\bf{\mathcal{C}}}=\{C_{1},\ldots,C_{n}\},$ where
\[C_{1}=\ldots=C_{n-i_{2}-1}=[n]\setminus \sigma^{t_{2}}[i_{2}],\]
\[C_{n-i_{2}}=\ldots =C_{i_{1}}=C_{n}=[n], \ \ \ \ \ \]
\[C_{i_{1}+1}=\ldots =C_{n-1}=[n]\setminus [i_{1}]. \ \ \ \ \]
The polymatroidal diagram associated is the following:

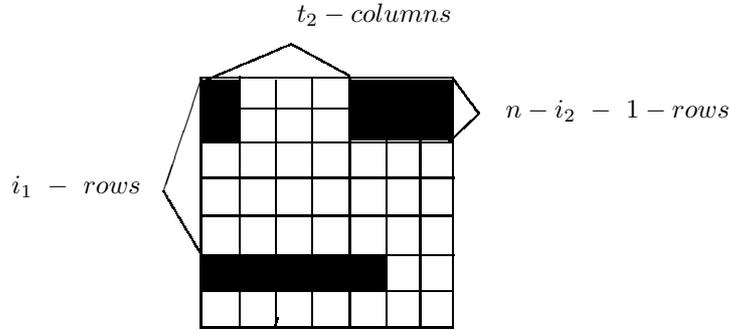
\unitlength 1mm 
\linethickness{0.4pt}
\ifx\plotpoint\undefined\newsavebox{\plotpoint}\fi 
\begin{picture}(60,70)(0,0)
\put(45,22.25){\framebox(33,33)[cc]{}}
\put(50,55.5){\line(0,-1){33.25}}
\put(54.75,55){\line(0,-1){32.25}}
\put(59.5,55.25){\line(0,-1){32.25}}
\multiput(55,23.5)(-.033333,-.1){15}{\line(0,-1){.1}}
\put(59.5,23.5){\line(0,-1){1.5}}
\put(64.5,55.25){\line(0,-1){32.75}}
\put(69.25,55){\line(0,-1){32.75}}
\put(73.75,55){\line(0,-1){32.5}}
\put(45,51.25){\line(1,0){32.75}}
\put(44.75,46.75){\line(1,0){33.25}}
\put(45,42){\line(1,0){33.25}}
\put(45.25,37){\line(1,0){32.75}}
\put(45.25,31.75){\line(1,0){33}}
\put(45,27){\line(1,0){33}}
\put(45.25,27.25){\rule{24\unitlength}{4.25\unitlength}}
\put(45.25,46.75){\rule{4.75\unitlength}{8.25\unitlength}}
\put(64.75,47){\rule{13.25\unitlength}{8\unitlength}}
\put(45,55.25){\line(-1,-3){5}}
\multiput(40,40.25)(.03355705,-.05704698){149}{\line(0,-1){.05704698}}
\multiput(45,54.75)(.07885906,.03355705){149}{\line(1,0){.07885906}}
\multiput(56.75,59.75)(.06150794,-.03373016){126}{\line(1,0){.06150794}}
\multiput(78,55.25)(.03365385,-.04567308){104}{\line(0,-1){.04567308}}
\multiput(81.5,50.5)(-.03571429,-.03348214){112}{\line(-1,0){.03571429}}
\put(57.5,62.75){$t_{2}-columns$}
\put(85,50){$n-i_{2} \ - \ 1-rows$}
\put(20,40){$i_{1} \ - \ rows$}
\end{picture}

It is easy to see that the base ring associated to the transversal polymatroid
${\bf{\mathcal{P}}}$ presented by ${\bf{\mathcal{C}}}$
is generated by the following set of monomials:
\[\{x_{i_{2}+t_{2}-n+1},\ldots,x_{t_{2}}\}^{i_{1}+1-k}
\{x_{1},\ldots, x_{i_{2}+t_{2}-n}, x_{t_{2}+1},\ldots,x_{i_{1}}\}^{k-s}
\{x_{i_{1}+1},\ldots,x_{n}\}^{n-1-i_{1}+ s}\]
for any $0\leq k \leq i_{1}+i_{2}-n+2$ and $0\leq s \leq k.$
Since $i_{2}+t_{2}\geq n$ and $0\leq s \leq k \leq i_{1}+i_{2}-n+2$
it follows that for any $x^{\alpha}\in K[{\bf{\mathcal{C}}}]$
we have $\alpha_{1}+\ldots+\alpha_{i_{1}}=i_{1}+1-s\leq i_{1}+1$,
$\alpha_{1}+\ldots +\alpha_{i_{2}+t_{2}-n}+\alpha_{t_{2}+1}+
\ldots+\alpha_{n}=n-1-i_{1}+k \leq i_{2}+1$
and thus, $K[{\bf{\mathcal{C}}}]\subset K[A\cap B].$\\
Conversely, \ if $\alpha \in A\cap B$ then $\alpha_{1}+\ldots+
\alpha_{i_{1}}\leq i_{1}+1$,
$\alpha_{1}+\ldots +\alpha_{i_{2}+t_{2}-n}+\alpha_{t_{2}+1}+\ldots+
\alpha_{n}\leq i_{2}+1$
and thus there exists $0\leq k \leq i_{1}+i_{2}-n+2$ and
$0\leq s \leq k$ such that \[x^{\alpha}\in \{x_{i_{2}+t_{2}-n+1},
\ldots,x_{t_{2}}\}^{i_{1}+1-k}
\{x_{1},\ldots, x_{i_{2}+t_{2}-n}, x_{t_{2}+1},\ldots,x_{i_{1}}\}^{k-s}
\{x_{i_{1}+1},\ldots,x_{n}\}^{n-1-i_{1}+ s}.\]
Thus, $K[{\bf{\mathcal{C}}}]\supset K[A\cap B]$ and so
$K[{\bf{\mathcal{C}}}]= K[A\cap B].$

$"\Rightarrow"$ \ Now suppose that there exists a transversal
polymatroid ${\bf{\mathcal{P}}}$ given by
${\bf{\mathcal{C}}}=\{C_{1},\ldots,C_{n}\}$ such that its associated
base ring is $K[A\cap B]$ and we will prove that
$i_{2}\in \{1,\ldots,i_{1}-t_{2}\}\cup \{n-t_{2},\ldots,n-2\}.$\\
Suppose, on the contrary, let $i_{1}+1-t_{2}\leq i_{2}\leq n-t_{2}-1.$
We have two cases to study:\\
${\bf Case} \ 1'.$ If $n-i_{1}-1\leq i_{2}+1,$ then since
$(i_{1}+1)e_{1}+(n-i_{1}-1)e_{k}\in {\bf{\mathcal{P}}}$
for any $i_{1}+1\leq k \leq n$ and
$(i_{1}+1)e_{1}+ e_{s}+(n-i_{1}-2)e_{k}\notin {\bf{\mathcal{P}}}$
for any $2\leq s \leq i_{1}$ and $i_{1}+1\leq k \leq n,$
we may assume $1\in C_{1},\ldots, 1\in C_{i_{1}}, 1\in C_{n}$
and $C_{i_{1}+1}=\ldots = C_{n-1}=[n]\setminus[i_{1}].$
If $i_{1}\leq i_{2},$ then since $(i_{1}+1)e_{t_{2}+1}+
(n-i_{1}-1)e_{t_{2}+i_{2}+1}\in {\bf{\mathcal{P}}},$
we may assume $t_{2}+1\in C_{1},\ldots, t_{2}+1\in C_{i_{1}},
t_{2}+1\in C_{n}.$ Then $(i_{1}+1)e_{t_{2}+1}+(n-i_{1}-1)
e_{i_{1}+1}\in {\bf{\mathcal{P}}},$ which is false.
If $i_{1}>i_{2},$ then since $(i_{2}+1)e_{t_{2}+1}+(n-i_{1}-1)
e_{t_{2}+i_{2}+1}\in {\bf{\mathcal{P}}},$ we may assume
$t_{2}+1\in C_{1},\ldots, t_{2}+1\in C_{i_{2}}, t_{2}+1\in C_{n}.$
Then $(i_{2}+1)e_{t_{2}+1}+(i_{1}-i_{2})e_{1}+(n-i_{1}-1)
e_{i_{1}+1}\in {\bf{\mathcal{P}}},$ which is false.\\
${\bf Case} \ 2'.$ If $n-i_{1}-1>i_{2}+1,$ then since $(i_{1}+1)e_{1}+
(i_{2}+1)e_{i_{1}+1}+(n-i_{1}-i_{2}-2)e_{k}\in {\bf{\mathcal{P}}}$
for any $t_{2}+i_{2}+1\leq k\leq n$ and $(i_{1}+1)e_{1}+
e_{s}+(i_{2}+1)e_{i_{1}+1}+(n-i_{1}-i_{2}-3)e_{k}\notin
{\bf{\mathcal{P}}}$ for any
$1\leq s \leq i_{1}$ and $t_{2}+i_{2}+1\leq k \leq n,$
we may assume $1\in C_{1},\ldots, 1\in C_{i_{1}}, 1\in C_{n},$
$C_{i_{1}+1}=\ldots = C_{i_{1}+i_{2}+1}=[n]\setminus[i_{1}]$ and
$C_{i_{1}+i_{2}+2}=\ldots=C_{n-1}=[n]\setminus[t_{2}+i_{2}].$
If $i_{1}\leq i_{2},$ then since
$(i_{1}+1)e_{t_{2}+1}+(n-i_{1}-1)e_{t_{2}+i_{2}+1}\in
{\bf{\mathcal{P}}},$ we may assume $t_{2}+1\in C_{1},\ldots,
t_{2}+1\in C_{i_{1}}, t_{2}+1\in C_{n}.$
Then $(i_{1}+1)e_{t_{2}+1}+(i_{2}+1)e_{i_{1}+1}+
(n-i_{1}-i_{2}-2)e_{t_{2}+i_{2}+1}\in
{\bf{\mathcal{P}}},$ which is false. If $i_{1}>i_{2},$
then since $(i_{2}+1)e_{t_{2}+1}+(i_{1}-i_{2})e_{1}+
(n-i_{1}-1)e_{t_{2}+i_{2}+1}\in {\bf{\mathcal{P}}},$
we may assume $t_{2}+1\in C_{1},\ldots,
t_{2}+1\in C_{i_{2}}, t_{2}+1\in C_{n}.$ Then
$(i_{1}-i_{2})e_{1}+(i_{2}+1)e_{t_{2}+1}+(i_{2}+1)e_{t_{2}+
i_{2}}+(n-i_{1}-i_{2}-2)e_{t_{2}+i_{2}+1}\in
{\bf{\mathcal{P}}},$ which is false.
\end{proof}

\begin{lemma}
Let $A$ and $B$ be like above. If $i_{1}\geq 2,$
$i_{1}+1\leq t_{2} \leq n-1,$  then  the
$K-algebra$ $K[A\cap B]$ is the base ring associated to some
transversal polymatroid if and only if
$i_{2}\in \{1,\ldots,n-t_{2}\}
\cup \{n-t_{2}+i_{1},\ldots,n-2\}.$
\end{lemma}
\begin{proof}
$"\Leftarrow"$ \
Let $i_{2}\in \{1,\ldots,n-t_{2}\}
\cup \{n-t_{2}+i_{1},\ldots,n-2\}.$ We will prove that
there exists a transversal polymatroid ${\bf{\mathcal{P}}}$
presented by  ${\bf{\mathcal{C}}}=\{C_{1},\ldots,C_{n}\}$
such that its associated base ring is $K[A\cap B].$
We distinct three cases to study:\\
${\bf Case} \ 1.$ If $i_{2}+t_{2}\leq n$ and $i_{1}+1+i_{2}\neq n,$
then let ${\bf{\mathcal{P}}}$ be the transversal polymatroid
presented by  ${\bf{\mathcal{C}}}=\{C_{1},\ldots,C_{n}\},$ where
\[C_{1}=\ldots =C_{i_{1}}=C_{n}=[n]\setminus \sigma^{t_{2}}[i_{2}],\]
\[C_{i_{1}+1}=\ldots =C_{i_{1}+i_{2}+1}=[n]\setminus [i_{1}],\ \]
\[\ \ \ \ \ \ \ \ \ \ \ \ C_{i_{1}+i_{2}+2}=\ldots =C_{n-1}=[n]
\setminus ([i_{1}]\cup \sigma^{t_{2}}[i_{2}]).\]
The polymatroidal diagram associated is the following:

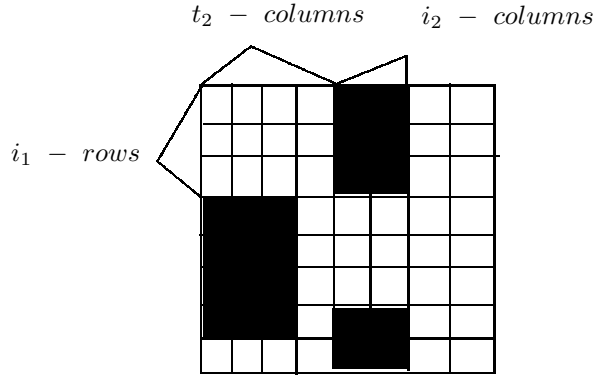
\unitlength 1mm 
\linethickness{0.4pt}
\ifx\plotpoint\undefined\newsavebox{\plotpoint}\fi 
\begin{picture}(93,65.75)(0,10)
\put(53.75,19.25){\framebox(38.5,38)[cc]{}}
\put(57.75,57.5){\line(0,-1){38.25}}
\put(61.75,57){\line(0,-1){37.75}}
\put(66.25,57){\line(0,-1){38}}
\put(71.25,57){\line(0,-1){37.25}}
\put(76,57.25){\line(0,-1){37.5}}
\put(81,57.25){\line(0,-1){37.75}}
\put(86.5,57.5){\line(0,-1){38.25}}
\put(54,52.25){\line(1,0){38.25}}
\put(53.75,48){\line(1,0){39.25}}
\put(53.5,42.5){\line(1,0){38.75}}
\put(53.5,37.5){\line(1,0){38.75}}
\put(53.75,33.25){\line(1,0){38.5}}
\put(53.75,28.25){\line(1,0){38.25}}
\put(53.75,23.75){\line(1,0){38.25}}
\put(54,23.75){\rule{12\unitlength}{18.5\unitlength}}
\put(71.25,43){\rule{9.75\unitlength}{14\unitlength}}
\put(71,19.75){\rule{9.75\unitlength}{8\unitlength}}
\multiput(53.75,42.5)(-.04078014,.03368794){141}{\line(-1,0){.04078014}}
\multiput(48,47.25)(.03370787,.05758427){178}{\line(0,1){.05758427}}
\multiput(53.5,57.25)(.04326923,.03365385){156}{\line(1,0){.04326923}}
\multiput(60.25,62.5)(.07550336,-.03355705){149}{\line(1,0){.07550336}}
\multiput(71.5,57.5)(.08258929,.03348214){112}{\line(1,0){.08258929}}
\put(80.75,61.25){\line(0,-1){4}}
\put(28.75,47.25){$i_{1} \ - \ rows$}
\put(52.5,65.75){$t_{2}\ - \ columns$}
\put(82.75,65.5){$i_{2} \ - \ columns$}
\end{picture}

It is easy to see that the base ring $K[{\bf{\mathcal{C}}}]$
associated to the transversal polymatroid ${\bf{\mathcal{P}}}$
presented by ${\bf{\mathcal{C}}}$ is generated by the
following set of monomials:
\[\{x_{1},\ldots,x_{i_{1}}\}^{i_{1}+1-k}
\{x_{t_{2}+1},\ldots,x_{t_{2}+i_{2}}\}^{i_{2}+1-s}
\{x_{i_{1}+1},\ldots,x_{t_{2}},x_{t_{2}+i_{2}+1},\ldots,x_{n}\}
^{n-i_{1}-i_{2}-2+k+s}\] for any $0\leq k\leq i_{1}+1$ and
$0\leq s\leq i_{2}+1.$
If $x^{\alpha}\in K[{\bf{\mathcal{C}}}],$
$\alpha=(\alpha_{1},\ldots,\alpha_{n})\in {\mathbb{N}}^{n},$
then there exists $0\leq k \leq i_{1}+1$ and $0\leq s \leq i_{2}+1$
such that
\[\alpha_{t_{2}+1}+\ldots + \alpha_{t_{2}+i_{2}}= i_{2}+1-s
\ , \ \ \alpha_{1}+\ldots + \alpha_{i_{1}}= i_{1}+1-k\]
and thus, $K[{\bf{\mathcal{C}}}]\subset K[A\cap B].$
Conversely, \ if $\alpha \in A\cap B$ then
$\alpha_{t_{2}+1}+\ldots + \alpha_{t_{2}+i_{2}}\leq i_{2}+1
\ , \ \ \alpha_{1}+\ldots + \alpha_{i_{1}}\leq i_{1}+1;$
thus there exists $0\leq k \leq i_{1}+1$ and $0\leq s
\leq i_{2}+1$ such that
\[\alpha_{t_{2}+1}+\ldots + \alpha_{t_{2}+i_{2}}= i_{2}+1-s
\ , \ \ \alpha_{1}+\ldots + \alpha_{i_{1}}= i_{1}+1-k\] and
since $ | \ \alpha \ | = n$ it follows that
 \[x^{\alpha}\in \{x_{1},\ldots,
x_{i_{1}}\}^{i_{1}+1-k}
\{x_{t_{2}+1},\ldots,x_{t_{2}+i_{2}}\}^{i_{2}+1-s}
\{x_{i_{1}+1},\ldots,x_{t_{2}},x_{t_{2}+i_{2}+1},\ldots,x_{n}\}
^{n-i_{1}-i_{2}-2+k+s}.\]
Thus, $K[{\bf{\mathcal{C}}}]\supset K[A\cap B]$ and so
$K[{\bf{\mathcal{C}}}]= K[A\cap B].$\\
${\bf Case} \ 2.$ If $i_{2}+t_{2}\leq n$ and $i_{1}+1+i_{2}=n,$ then
$t_{2}=i_{1}+1$ and let ${\bf{\mathcal{P}}}$
be the transversal polymatroid presented by
${\bf{\mathcal{C}}}=\{C_{1},\ldots,C_{n}\},$ where
\[C_{1}=\ldots =C_{i_{1}}=[n]\setminus \sigma^{t_{2}}[i_{2}],\]
\[\ \  C_{i_{1}+1}= \ldots =C_{n-1}=[n]\setminus [i_{1}],\]
\[\ \ \  C_{n}=[n].\ \ \ \ \ \ \ \ \ \ \ \ \ \ \ \ \ \ \ \ \ \
\ \ \ \ \ \ \ \  \]
The polymatroidal diagram associated is the following:

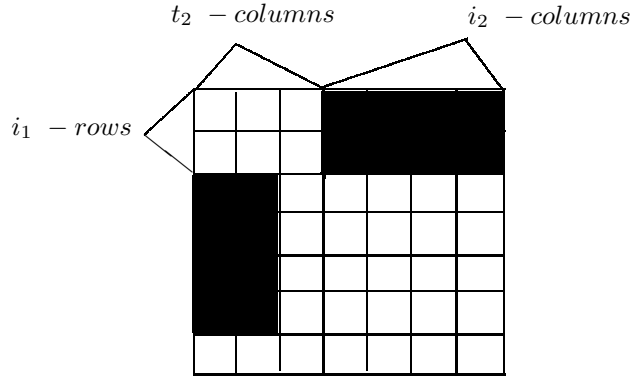
\unitlength 1mm 
\linethickness{0.4pt}
\ifx\plotpoint\undefined\newsavebox{\plotpoint}\fi 
\begin{picture}(0,70)(10,10)
\put(54,18.25){\framebox(40.75,37.75)[cc]{}}
\put(59.5,55.5){\line(0,-1){37}}
\put(65.25,56){\line(0,-1){37.25}}
\put(71,56){\line(0,-1){37.5}}
\put(76.75,55.75){\line(0,-1){37}}
\put(82.5,55.75){\line(0,-1){37.5}}
\put(88.5,56){\line(0,-1){37.75}}
\put(54,50.5){\line(1,0){41}}
\put(54.25,44.75){\line(1,0){40.5}}
\put(54,39.75){\line(1,0){40.75}}
\put(54,34){\line(1,0){40.75}}
\put(53.75,29.25){\line(1,0){41}}
\put(54,23.5){\line(1,0){40.5}}
\put(53.75,23.75){\rule{.5\unitlength}{1\unitlength}}
\multiput(70.75,45)(.03125,-.03125){8}{\line(0,-1){.03125}}
\put(74.25,45){\rule{.25\unitlength}{.25\unitlength}}
\put(53.75,23.75){\rule{11.25\unitlength}{20.75\unitlength}}
\multiput(54,56)(-.03651685,-.03370787){178}{\line(-1,0){.03651685}}
\put(47.5,50){\line(5,-4){6.25}}
\multiput(54,55.75)(.03353659,.03810976){164}{\line(0,1){.03810976}}
\multiput(59.5,62)(.06578947,-.03362573){171}{\line(1,0){.06578947}}
\multiput(70.75,56.25)(.09974093,.03367876){193}{\line(1,0){.09974093}}
\put(30,50){$i_{1} \ -rows$}
\put(51,65){$t_{2} \ -columns$}
\put(90,64.75){$i_{2} \ -columns$}
\multiput(89.75,62.5)(.03355705,-.04530201){149}{\line(0,-1){.04530201}}
\put(71,44.25){\rule{.25\unitlength}{.75\unitlength}}
\put(70.75,45){\rule{23.75\unitlength}{10.75\unitlength}}
\end{picture}

It is easy to see that the base ring $K[{\bf{\mathcal{C}}}]$
associated to the transversal polymatroid ${\bf{\mathcal{P}}}$
presented by ${\bf{\mathcal{C}}}$ is generated by the
following set of monomials:
\[\{x_{1},\ldots,x_{i_{1}}\}^{i_{1}+1-k}x_{i_{1}+1}^{n-i_{1}-1+k-s}
\{x_{i_{1}+2},\ldots,x_{n}\}^{s}\] for any
$0\leq k\leq i_{1}+1$ and $0\leq s\leq n-i_{1}.$
If $x^{\alpha}\in K[{\bf{\mathcal{C}}}],$ $\alpha=
(\alpha_{1},\ldots,\alpha_{n})\in {\mathbb{N}}^{n},$
then there exists $0\leq k \leq i_{1}+1$ and
$0\leq s \leq i_{2}+1(=n-i_{1})$ such that
\[\alpha_{t_{2}+1}+\ldots + \alpha_{t_{2}+i_{2}}=
\alpha_{i_{1}+2}+\ldots+x_{n}=s\leq i_{2}+1 \ , \  \
\alpha_{1}+\ldots + \alpha_{i_{1}}= i_{1}+1-k\] and thus,
$K[{\bf{\mathcal{C}}}]\subset K[A\cap B].$
Conversely, \ if $\alpha \in A\cap B$ then
$\alpha_{t_{2}+1}+\ldots + \alpha_{t_{2}+i_{2}}=
\alpha_{i_{1}+2}+\ldots+x_{n}\leq i_{2}+1 \ and \
\alpha_{1}+\ldots + \alpha_{i_{1}}\leq i_{1}+1;$
thus there exists $0\leq k \leq i_{1}+1$ and
$0\leq s \leq i_{2}+1$ such that \[x^{\alpha}\in
\{x_{1},\ldots,x_{i_{1}}\}^{i_{1}+1-k}x_{i_{1}+1}^
{n-i_{1}-1+k-s}\{x_{i_{1}+2},\ldots,x_{n}\}^{s}.\]
Thus, $K[{\bf{\mathcal{C}}}]\supset K[A\cap B]$ and
so $K[{\bf{\mathcal{C}}}]= K[A\cap B].$\\
${\bf Case} \ 3.$ If $i_{2}+t_{2}>n,$ then let ${\bf{\mathcal{P}}}$
be the transversal polymatroid presented by
${\bf{\mathcal{C}}}=\{C_{1},\ldots,C_{n}\},$ where
\[C_{1}=\ldots =C_{i_{1}}=C_{n}=[n],\]
\[\ \ \ \ \ \ \ \ \ \ \ \ \ \ \ \ \ C_{i_{1}+1}=
\ldots =C_{i_{1}+n-i_{2}-1}=[n]\setminus \sigma^{t_{2}}[i_{2}],\]
\[\ \ \ \ \ \ \ \ \ C_{i_{1}+n-i_{2}}=\ldots=C_{n-1}=[n]\setminus[i_{1}].\]
The polymatroidal diagram associated is the following:

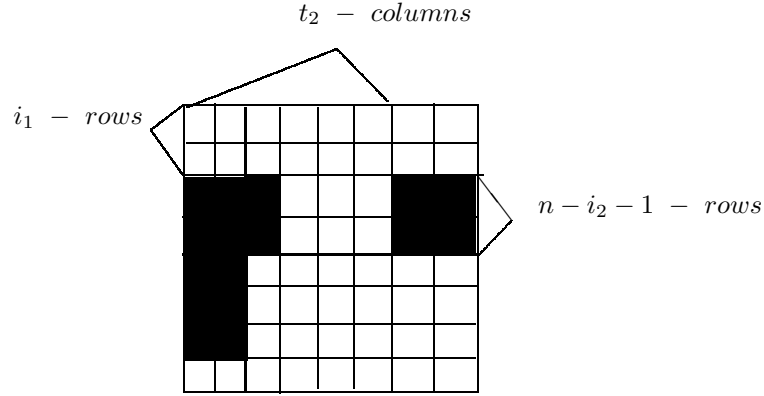
\unitlength 1mm 
\linethickness{0.4pt}
\ifx\plotpoint\undefined\newsavebox{\plotpoint}\fi 
\begin{picture}(0,30)(-20,45)
\put(32.5,14.25){\framebox(38.5,38)[cc]{}}
\put(36.5,52.5){\line(0,-1){38.25}}
\put(40.5,52){\line(0,-1){37.75}}
\put(45,52){\line(0,-1){38}}
\put(50,52){\line(0,-1){37.25}}
\put(54.75,52.25){\line(0,-1){37.5}}
\put(59.75,52.25){\line(0,-1){37.75}}
\put(65.25,52.5){\line(0,-1){38.25}}
\put(32.75,47.25){\line(1,0){38.25}}
\put(32.5,43){\line(1,0){39.25}}
\put(32.25,37.5){\line(1,0){38.75}}
\put(32.25,32.5){\line(1,0){38.75}}
\put(32.5,28.25){\line(1,0){38.5}}
\put(32.5,23.25){\line(1,0){38.25}}
\put(32.5,18.75){\line(1,0){38.25}}
\put(32.5,18.5){\rule{8.25\unitlength}{24.25\unitlength}}
\put(40.5,32.75){\rule{4.25\unitlength}{10\unitlength}}
\put(59.75,32.5){\rule{11\unitlength}{10.25\unitlength}}
\multiput(32.25,43)(-.03373016,.04761905){126}{\line(0,1){.04761905}}
\multiput(28,49)(.04381443,.03350515){97}{\line(1,0){.04381443}}
\multiput(32.75,52)(.08586957,.03369565){230}{\line(1,0){.08586957}}
\multiput(52.5,59.75)(.03358209,-.03482587){201}{\line(0,-1){.03482587}}
\put(71,43){\line(3,-4){4.5}}
\multiput(75.5,37)(-.03358209,-.03544776){134}{\line(0,-1){.03544776}}
\put(47.5,63.5){$t_{2} \ - \ columns$}
\put(10,49.75){$i_{1} \ - \ rows$}
\put(79,38){$n-i_{2}-1 \ - \ rows$}
\end{picture}
\[\] \[\] \[\] \[\] \[\] \[\] 
Since $i_{2}+t_{2}>n$ it follows that $i_{2}+t_{2}\geq n+i_{1}$
and so $i_{2}-i_{1}\geq n-t_{2}\geq 1.$
It is easy to see that the base ring $K[{\bf{\mathcal{C}}}]$
associated to the transversal polymatroid ${\bf{\mathcal{P}}}$
presented by ${\bf{\mathcal{C}}}$ is generated by the
following set of monomials:
\[\{x_{1},\ldots,x_{i_{1}}\}^{i_{1}+1-k}
\{x_{i_{1}+1},\ldots,x_{i_{2}+t_{2}-n},x_{t_{2}+1},\ldots,x_{n}\}
^{i_{2}-i_{1}+k-s}
\{x_{i_{2}+t_{2}-n+1},\ldots,x_{t_{2}}\}^{n-i_{2}-1+s}\] for any
$0\leq k \leq i_{1}+1$ and $0\leq s \leq i_{2}-i_{1}+k.$ If
$x^{\alpha}\in K[{\bf{\mathcal{C}}}],$
$\alpha=(\alpha_{1},\ldots,\alpha_{n})\in {\mathbb{N}}^{n},$ then
there exists $0\leq k \leq i_{1}+1$ and $0\leq s \leq i_{2}-i_{1}+k$
such that
\[\alpha_{1}+\ldots + \alpha_{i_{1}}= i_{1}+1-k \ , \
\ \alpha_{1}+\ldots+\alpha_{i_{2}+t_{2}-n}+
\alpha_{t_{2}+1}+\ldots+\alpha_{n}= i_{2}+1-s\] and thus,
$K[{\bf{\mathcal{C}}}]\subset K[A\cap B].$ Conversely, \ if $\alpha
\in A\cap B$ then $\alpha_{1}+ \ldots + \alpha_{i_{1}}\leq i_{1}+1 \
and \ \alpha_{1}+ \ldots+\alpha_{i_{2}+t_{2}-n}+\alpha_{t_{2}+1}+
\ldots+\alpha_{n}\leq i_{2}+1.$ Then there exists $0\leq k \leq
i_{1}+1$ and $0\leq s \leq i_{2}-i_{1}+k$ such that
\[\alpha_{1}+\ldots + \alpha_{i_{1}}= i_{1}+1-k \ , \ \
\alpha_{1}+\ldots+\alpha_{i_{2}+t_{2}-n}+
\alpha_{t_{2}+1}+\ldots+\alpha_{n}= i_{2}+1-s\] and since $|\ \alpha
\ |=n$ it follows that
 \[x^{\alpha}\in \{x_{1},\ldots,x_{i_{1}}\}^{i_{1}+1-k}
\{x_{i_{1}+1},\ldots,x_{i_{2}+t_{2}-n},x_{t_{2}+1},\ldots,x_{n}\}
^{i_{2}-i_{1}+k-s}
\{x_{i_{2}+t_{2}-n+1},\ldots,x_{t_{2}}\}^{n-i_{2}-1+s}.\]
Thus, $K[{\bf{\mathcal{C}}}]\supset K[A\cap B]$
and so $K[{\bf{\mathcal{C}}}]= K[A\cap B].$

$"\Rightarrow"$ \ Now suppose that there exists a transversal
polymatroid ${\bf{\mathcal{P}}}$ presented by
${\bf{\mathcal{C}}}=\{C_{1},\ldots,C_{n}\}$ such that its associated
base ring  is $K[A\cap B]$ and we will prove that $i_{2}\in
\{1,\ldots,n-t_{2}\} \cup \{n-t_{2}+i_{1},\ldots,n-2\}.$ Suppose, on
the contrary, let $i_{2}\in \{n-t_{2}+1,\ldots,n-t_{2}+i_{1}-1\}.$
We may assume that $1\notin C_{i_{2}+t_{2}-n+1}, \ldots, 1\notin
C_{i_{1}}, 1\notin C_{i_{1}+1}, \ldots, 1\notin C_{n-1}.$ If $i_{1}<
i_{2},$ then $x_{1}^{i_{1}+1}x_{i_{1}+1}^{n-i_{1}-1}\in K[A\cap B]$.
But $x_{1}^{i_{1}+1}x_{i_{1}+1}^{n-i_{1}-1} \notin
K[{\bf{\mathcal{C}}}]$ because the maximal power of $x_1$ in a
minimal generator of $K[C]$ is $\leq i_2+t_2-n+1\leq i_1$, which is
false. If $i_{1}\geq i_{2},$ then
$x_{1}^{i_{2}+1}x_{i_{1}+1}^{n-i_{2}-1}\in K[A\cap B]$. But
$x_{1}^{i_{2}+1}x_{i_{1}+1}^{n-i_{2}-1} \notin
K[{\bf{\mathcal{C}}}]$ because the maximal power of $x_1$ in a
minimal generator of $K[C]$ is $\leq i_2+t_2-n+1\leq i_2$, which is
false. Thus, $i_{2}\in \{1,\ldots,n-t_{2}\} \cup
\{n-t_{2}+i_{1},\ldots,n-2\}.$
\end{proof}

\begin{lemma}
Let $A$ and $B$ like above. If either $i_{1}=1$ and  $0\leq t_{2}
\leq n-1,$ or $i_{1}\geq 2$ and $t_{2}=i_{1},$ or $i_{1}\geq 2$ and
$t_{2}=0,$   then  the $K-algebra$ $K[A\cap B]$ is the base ring
associated to some transversal polymatroid.
\end{lemma}
\begin{proof}
We have three cases to study:\\
${\bf Case} \ 1.$ \ $i_{1}=1$ and  $0\leq t_{2} \leq n-1.$ \
Then we distinct five subcases:\\
${\bf Subcase} \ 1.a.$ If $t_{2}=0,$ then we find a transversal polymatroid
${\bf{\mathcal{P}}}$
like in the Subcase $3.b.$ when $i_{1}=1.$\\
${\bf Subcase} \ 1.b.$ If $t_{2}>0 \ and \ t_{2}+i_{2} \leq n \ with \
i_{2}\neq n-2, \ n-3,$ then
let ${\bf{\mathcal{P}}}$ be the transversal polymatroid
presented by  ${\bf{\mathcal{C}}}=\{C_{1},\ldots,C_{n}\},$
where \[C_{1}=C_{n}=[n]\setminus \sigma^{t_{2}}[i_{2}],
\ \ \ \ \ \ \ \ \ \ \ \ \ \ \ \ \ \ \ \ \ \ \ \ \]
\[C_{2}=\ldots =C_{i_{2}+2}=[n]\setminus[1], \ \ \ \ \ \
\ \ \ \ \ \ \ \ \ \ \ \ \]
\[ \ \ C_{i_{2}+3}=\ldots=C_{n-1}=[n]\setminus \{\{1\}\cup
\sigma^{t_{2}}[i_{2}]\}.\]
It is easy to see that the polymatroid ${\bf{\mathcal{P}}}$
is the same like in $Lemma \ 4.3$
when $i_{2}+t_{2}\leq n \ and \ i_{1}+1+i_{2}\neq n.$
Thus $K[A\cap B]=K[\bf{\mathcal{C}}].$\\
${\bf Subcase} \ 1.c.$ If $t_{2}>0 \ and \ t_{2}+i_{2} \leq n \
with \ i_{2}= n-2,$ then
let ${\bf{\mathcal{P}}}$ be the transversal polymatroid
presented by
${\bf{\mathcal{C}}}=\{C_{1},\ldots,C_{n}\},$ where
\[ \ \ \ C_{1}=[n]\setminus \sigma^{t_{2}}[n-2]\ , \ C_{n}=[n],
\ \ \ \ \ \ \ \ \ \ \ \ \ \ \ \]
\[C_{2}=\ldots =C_{n-1}=[n]\setminus[1]. \ \ \ \ \ \ \ \
\ \ \ \ \ \ \ \ \ \]
It is easy to see that the polymatroid ${\bf{\mathcal{P}}}$
is the same like in $Lemma \ 4.3.$
when $i_{2}+t_{2}\leq n \ and \ i_{1}+1+i_{2} = n.$
Thus $K[A\cap B]=K[\bf{\mathcal{C}}].$\\
${\bf Subcase} \ 1.d.$ If $t_{2}>0 \ and \ t_{2}+i_{2} \leq n \ with
\ i_{2}= n-3,$ then let ${\bf{\mathcal{P}}}$ be the transversal
polymatroid presented by
${\bf{\mathcal{C}}}=\{C_{1},\ldots,C_{n}\},$ where \[\
C_{1}=C_{n}=[n]\setminus \sigma^{t_2}[n-3], \ \ \ \ \ \ \ \ \ \ \ \
\ \ \ \ \ \ \ \ \]
\[\ \ \ \ C_{2}=\ldots =C_{n-1}=[n]\setminus[1], \ \ \
\ \ \ \ \ \ \ \ \ \ \ \ \ \ \ \ \ .\] It is easy to see that the polymatroid
${\bf{\mathcal{P}}}$ is the same like in $Lemma \ 4.3.$ when
$i_{2}+t_{2}\leq n \ and \ i_{1}+1+i_{2}\not = n.$
Thus $K[A\cap B]=K[\bf{\mathcal{C}}].$\\
${\bf Subcase} \ 1.e.$ If $t_{2}>0$ and $t_{2}+i_{2}>n,$ then
let ${\bf{\mathcal{P}}}$ be the transversal polymatroid
presented by \\
${\bf{\mathcal{C}}}=\{C_{1},\ldots,C_{n}\},$ where
\[ \ \ \ C_{1}=C_{n}=[n], \ \ \ \ \ \ \ \ \ \ \ \ \ \ \ \
 \ \ \ \ \ \ \ \ \ \ \ \ \ \ \ \ \ \ \ \  \]
\[ \ \ \ \ C_{2}=\ldots =C_{n-i_{2}}=[n]\setminus \sigma^{t_{2}}
[i_{2}], \ \ \ \ \ \ \ \ \ \ \ \ \ \ \ \]
\[ \ \ \ C_{n-i_{2}+1}=\ldots=C_{n-1}=[n]\setminus \{1\}.
\ \ \ \ \ \ \ \ \ \ \ \]
It is easy to see that the polymatroid ${\bf{\mathcal{P}}}$
is the same like in $Lemma \ 4.3$
when $i_{2}+t_{2}> n \ and \ i_{1}=1.$ Thus $K[A\cap B]=K[\bf{\mathcal{C}}].$\\
${\bf Case} \ 2.$ \ $i_{1}\geq 2$ and  $t_{2}=i_{1}.$\
Then we distinct three subcases:\\
${\bf Subcase} \ 2.a.$ If $i_{2}+t_{2}<n-1,$ then
let ${\bf{\mathcal{P}}}$ be the transversal polymatroid presented by \\
${\bf{\mathcal{C}}}=\{C_{1},\ldots,C_{n}\},$ where
\[C_{1}=\ldots=C_{i_{1}}=C_{n}=[n]\setminus \sigma^{t_{2}}[i_{2}], 
\ \ \ \ \ \]
\[C_{i_{1}+1}=\ldots =C_{i_{1}+i_{2}+1}=[n]\setminus [i_{1}], 
\ \ \ \ \ \ \]
\[ \ \ C_{i_{1}+i_{2}+2}=\ldots=C_{n-1}=[n]\setminus[i_{1}+i_{2}]. 
\ \ \]
It is easy to see that the polymatroid ${\bf{\mathcal{P}}}$
is the same like in $Lemma \ 4.3$
when $i_{2}+t_{2}\leq n \ and \ i_{1}+1 +i_{2}\neq n.$
Thus $K[A\cap B]=K[\bf{\mathcal{C}}].$\\
${\bf Subcase} \ 2.b.$ If $i_{2}+t_{2}=n-1,$ then
let ${\bf{\mathcal{P}}}$ be the transversal polymatroid presented by
${\bf{\mathcal{C}}}=\{C_{1},\ldots,C_{n}\},$ where
\[C_{1}=\ldots=C_{i_{1}}=[n]\setminus \sigma^{t_{2}}[i_{2}], \
\ \ \ \ \ \ \ \ \ \ \ \]
\[ \ \ C_{i_{1}+1}=\ldots =C_{n-1}=[n]\setminus [i_{1}], 
\ \ \ \ \ \ \ \ \ \ \ \ \]
\[\ C_{n}=[n]. \ \ \ \ \ \ \ \ \ \ \ \ \ \ \ \ \ \ \ \ \ \ \ \ \ \
\ \ \ \ \ \ \ \ \ \ \ \ \ \]
It is easy to see that the polymatroid ${\bf{\mathcal{P}}}$
is the same like in $Lemma \ 4.3$
when $i_{2}+t_{2}\leq n \ and \ i_{1}+1 +i_{2}= n.$
Thus $K[A\cap B]=K[\bf{\mathcal{C}}].$\\
${\bf Subcase} \ 2.c.$ If $i_{2}+t_{2}\geq n,$ then
let ${\bf{\mathcal{P}}}$ be the transversal polymatroid presented by
${\bf{\mathcal{C}}}=\{C_{1},\ldots,C_{n}\},$ where
\[ \ \ \ \ \ \ \ \ \ \ \ C_{1}=\ldots=C_{n-i_{2}-1}=[n]\setminus
\sigma^{t_{2}}[i_{2}], \ \ \ \ \ \ \ \ \ \ \ \]
\[\ \ \ \ \ \ C_{n-i_{2}}=\ldots =C_{i_{1}}=C_{n}=[n], \ \ \ \ \ \ \
\ \ \ \ \]
\[  C_{i_{1}+1}=\ldots=C_{n-1}=[n]\setminus[i_{1}]. \ \ \ \ \]
It is easy to see that the polymatroid ${\bf{\mathcal{P}}}$
is the same like in $Lemma \ 4.2$
when $i_{2}+t_{2}>i_{1}.$ Thus $K[A\cap B]=K[\bf{\mathcal{C}}].$\\
${\bf Case} \ 3.$ \ $i_{1}\geq 2$ and  $t_{2}=0.$\ Then we distinct two subcases:\\
${\bf Subcase} \ 3.a.$ If $i_{2}\leq i_{1},$ then
let ${\bf{\mathcal{P}}}$ be the transversal polymatroid presented by
${\bf{\mathcal{C}}}=\{C_{1},\ldots,C_{n}\},$ where
\[C_{1}=\ldots=C_{i_{2}}=C_{n}=[n], \ \ \ \ \ \ \ \ \]
\[ \ \ \ \ C_{i_{2}+1}=\ldots =C_{i_{1}}=[n]\setminus [i_{2}], \ \
\ \ \ \ \ \ \ \]
\[ \ \ C_{i_{1}+1}=\ldots=C_{n-1}=[n]\setminus[i_{1}]. \ \ \ \ \ \]
It is easy to see that the polymatroid ${\bf{\mathcal{P}}}$ is
the same like in $Lemma \ 4.2$
when $i_{2}+t_{2}\leq i_{1}.$ Thus $K[A\cap B]=K[\bf{\mathcal{C}}].$\\
${\bf Subcase} \ 3.b.$ If $i_{2}> i_{1},$ then
let ${\bf{\mathcal{P}}}$ be the transversal polymatroid presented by
${\bf{\mathcal{C}}}=\{C_{1},\ldots,C_{n}\},$ where
\[C_{1}=\ldots=C_{i_{1}}=C_{n}=[n], \ \ \ \ \ \ \ \ \]
\[ \ \ \ \ C_{i_{1}+1}=\ldots =C_{i_{2}}=[n]\setminus [i_{1}], \ \ \
\ \ \ \ \ \ \]
\[ C_{i_{2}+1}=\ldots=C_{n-1}=[n]\setminus[i_{2}]. \ \ \]
The polymatroidal diagram associated is the following:

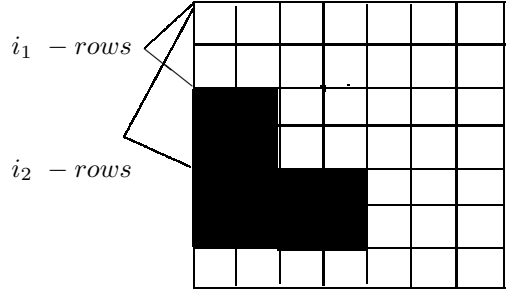
\unitlength 1mm 
\linethickness{0.4pt}
\ifx\plotpoint\undefined\newsavebox{\plotpoint}\fi 
\begin{picture}(0,70)(10,0)
\put(54,18.25){\framebox(40.75,37.75)[cc]{}}
\put(59.5,55.5){\line(0,-1){37}}
\put(65.25,56){\line(0,-1){37.25}}
\put(71,56){\line(0,-1){37.5}}
\put(76.75,55.75){\line(0,-1){37}}
\put(82.5,55.75){\line(0,-1){37.5}}
\put(88.5,56){\line(0,-1){37.75}}
\put(54,50.5){\line(1,0){41}}
\put(54.25,44.75){\line(1,0){40.5}}
\put(54,39.75){\line(1,0){40.75}}
\put(54,34){\line(1,0){40.75}}
\put(53.75,29.25){\line(1,0){41}}
\put(54,23.5){\line(1,0){40.5}}
\put(53.75,23.75){\rule{.5\unitlength}{1\unitlength}}
\multiput(70.75,45)(.03125,-.03125){8}{\line(0,-1){.03125}}
\put(74.25,45){\rule{.25\unitlength}{.25\unitlength}}
\put(53.75,23.75){\rule{11.25\unitlength}{20.75\unitlength}}
\multiput(54,56)(-.03651685,-.03370787){178}{\line(-1,0){.03651685}}
\put(47.5,50){\line(5,-4){6.25}}
\put(30,49.25){$i_{1} \ -rows$}
\put(71,44.25){\rule{.25\unitlength}{.75\unitlength}}
\multiput(54,55.5)(-.033636364,-.062727273){275}{\line(0,-1){.062727273}}
\multiput(44.75,38.25)(.07352941,-.03361345){119}{\line(1,0){.07352941}}
\put(30,33){$i_{2} \ - rows$}
\put(65,23.25){\rule{11.5\unitlength}{10.5\unitlength}}
\end{picture}

It is easy to see that the base ring $K[{\bf{\mathcal{C}}}]$
associated to the transversal polymatroid ${\bf{\mathcal{P}}}$
presented by ${\bf{\mathcal{C}}}$ is generated by the following
set of monomials:
\[\{x_{1},\ldots,x_{i_{1}}\}^{i_{1}+1-k}
\{x_{i_{1}+1},\ldots,x_{i_{2}}\}^{i_{2}-i_{1}+k-s}
\{x_{i_{2}+1},\ldots,x_{n}\}^{n-i_{2}+s-1}\] for any $0\leq k\leq
i_{1}+1$ and $0\leq s\leq i_{2}-i_{1}+k.$ If $x^{\alpha}\in
K[{\bf{\mathcal{C}}}],$ $\alpha=(\alpha_{1},\ldots,\alpha_{n})\in
{\mathbb{N}}^{n},$ then there exists $0\leq k \leq i_{1}+1$ and
$0\leq s \leq i_{2}-i_{1}+k$ such that
\[\alpha_{1}+\ldots + \alpha_{i_{2}}= i_{2}+1-s \ and \
\alpha_{1}+\ldots + \alpha_{i_{1}}= i_{1}+1-k\] and thus,
$K[{\bf{\mathcal{C}}}]\subset K[A\cap B].$ Conversely, \ if $\alpha
\in A\cap B$ then $\alpha_{1}+ \ldots + \alpha_{i_{2}}\leq i_{2}+1 \
and \ \alpha_{1}+\ldots + \alpha_{i_{1}}\leq i_{1}+1$ and so there
exists $0\leq k \leq i_{1}+1$ and $0\leq s \leq i_{2}-i_{1}+k$ such
that
\[\alpha_{1}+\ldots + \alpha_{i_{2}}= i_{2}+1-s \ and \
\alpha_{1}+\ldots + \alpha_{i_{1}}= i_{1}+1-k\] and since
$| \ \alpha \ |=n$ it follows that
\[\{x_{1},\ldots,x_{i_{1}}\}^{i_{1}+1-k}
\{x_{i_{1}+1},\ldots,x_{i_{2}}\}^{i_{2}-i_{1}+k-s}
\{x_{i_{2}+1},\ldots,x_{n}\}^{n-i_{2}+s-1}\] Thus,
$K[{\bf{\mathcal{C}}}]\supset K[A\cap B]$
and so $K[{\bf{\mathcal{C}}}]= K[A\cap B].$\\
\end{proof}

{\bf Acknowledgments:}
The author whould like to thank Professor Dorin Popescu for valuable suggestions
and comments during the preparation of this paper.

\vspace{2mm} \noindent {\footnotesize
\begin{minipage}[b]{8cm}
Alin \c{S}tefan, Assistant Professor\\
"Petroleum and Gas" University of Ploie\c{s}ti\\
Ploie\c{s}ti, Romania\\
E-mail:nastefan@upg-ploiesti.ro
\end{minipage}}

\begin{thebibliography}{99}


\bibitem{BH}  W. Bruns, J. Herzog, {\it Cohen-Macaulay rings,} Revised
Edition, Cambridge, 1997.

\bibitem{BK} W. Bruns, R. Koch, \textit{Normaliz}--a program for
 computing
normalizations of affine semigroups, 1998. Available via anonymous
ftp from: ftp.mathematik.Uni-Osnabrueck.DE/pub/osm/kommalg/software.

\bibitem{HH} J. Herzog, T. Hibi, Discrete polymatroids,  {\it J.
 Algebraic Combin.}, \textbf{16}(2002), no. 3, 239--268.

\bibitem{HHV} J. Herzog, T. Hibi, M. Vladoiu,  Ideals of fiber type
and polymatroids, Osaka J. Math. 42 (2005), 807-829.

\bibitem{SA} A. \c{S}tefan,  {\it A class of transversal polymatroids with
Gorenstein base ring,} Bull. Math. Soc. Sci. Math. Roumanie Tome 51(99)
No. 1, 2008.

\bibitem{V}  R. Villarreal, {\it Monomial Algebras,} Marcel Dekker,
 New-York, 2001.

\bibitem{V1} M. Vladoiu, Discrete polymatroids,  An. \c{S}t. Univ.
Ovidius, Constan\c ta, 14 (2006), 89-112.

\bibitem{V2} M. Vladoiu, Equidimensional and unmixed ideals of
Veronese type, to appear in Communications in Alg., arXiv:math.
AC/0611326.

\end{thebibliography}
\end{document}